\documentclass[11pt]{article}
\usepackage{amssymb,amsmath,amsthm}
\usepackage{epsfig}
\usepackage{verbatim}
\usepackage{graphicx}

\usepackage{color}

\usepackage{MnSymbol}

\newtheorem{thm}{Theorem}[section]

\newtheorem{prop}[thm]{Proposition}
\newtheorem{defi}[thm]{Definition}
\newtheorem{rem}[thm]{Remark}

\newcommand{\R}{\mathbb{R}}

\newcommand{\Z}{\Bbb{Z}}
\newcommand{\N}{\Bbb{N}}

\newcommand{\grad}{\nabla}

\newcommand{\eps}{\epsilon}

\newcommand{\heatinv}{(\partial_t-\Delta)_0^{-1}}

\begin{document}

\title{Global regularity of 2D density patches\\for inhomogeneous Navier-Stokes}

\author{Francisco Gancedo and Eduardo Garc\'ia-Ju\'arez}

\date{\today}

\maketitle

\begin{abstract}
This paper is about Lions' open problem on density patches \cite{LIONS}: whether inhomogeneous incompressible Navier-Stokes equations preserve the initial regularity of the free boundary given by density patches. Using classical Sobolev spaces for the velocity, we first establish the propagation of $C^{1+\gamma}$ regularity with $0<\gamma<1$ in the case of positive density. Furthermore, we go beyond to show the persistence of a geometrical quantity such as the curvature. In addition, we obtain a proof for $C^{2+\gamma}$ regularity. 
\end{abstract}


\setcounter{tocdepth}{1}

\section{Introduction}

We consider an incompressible inhomogeneous fluid in the whole space $\R^2$,
\begin{equation}\label{INH}
\begin{gathered}
\grad\cdot u=0,\\
\partial_t\rho+u\cdot\grad \rho=0,\\
\end{gathered}
\end{equation}
driven by Navier-Stokes equations
\begin{equation}
\label{NavierStokes}
\rho(\partial_t u+ u\cdot\grad u) =\Delta u-\grad p,
\end{equation}
where the unknowns $\rho, u, p$ represent the density, velocity field and pressure of the fluid.

In the case of positive density, the first results of existence of strong solutions for smooth initial data were proved by Ladyzhenskaya and Solonnikov \cite{LADY}. When $\rho\geq 0$ is allowed, Simon \cite{SIMON} proved the global existence of weak solutions with finite energy. Afterwards, this result was extended to the case with variable viscosity by Lions in \cite{LIONS}. There, the author proposed the so-called {\em{density patch problem}}: assuming $\rho_0=1_{D_0}$ for some domain $D_0\subset \R^2$, the question is whether or not $\rho(t)=1_{D(t)}$ for some domain $D(t)$ with the same regularity as the initial one. Theorem 2.1 in \cite{LIONS} ensures that the density remains as a patch preserving its volume, but gives no information about the persistence of regularity.

Previously to this problem, the analogous question in vortex patches in Euler equations arose great interest, due to the fact that several numerical results indicated the possible formation of finite time singularities. 
First Chemin \cite{CHEMIN} using paradifferential calculus and later Bertozzi and Constantin \cite{CONST} by a geometrical harmonic analysis approach finally solved the {\em{vortex patch problem}} proving the contrary: $C^{1+\gamma}$ vortex patches preserve their regularity in time.

On the other hand, the appearance of finite-time singularities has been proved in related scenarios. For the Muskat problem density patches have been shown to become singular in finite time \cite{FINITEWATERWAVES}, \cite{MUSKAT}. When vacuum is considered for Euler equations with gravity, \lq splash' singularities were shown in the so-called water wave problem \cite{FINITEEULER}. Later these results were extended to parabolic problems such as Muskat \cite{MUSKAT2} and Navier-Stokes \cite{FINITENAVIER}. Different proofs of these results can be found in \cite{COUTAND1}, \cite{COUTAND2}. In \cite{FEFFERMANABSENCE} it is shown that the presence of a second fluid precludes \lq splash' singularities in Euler equations with gravity and surface tension. 
See \cite{GANCEDOSTRAIN} for a different proof applied to the Muskat problem and also \cite{COUTAND3} for the result including vorticity in the bulk.

Global-in-time regularity has been extensively studied for Navier-Stokes free boundary problems considering the continuity of the stress tensor at the free boundary (see \cite{WEHAUSEN} and \cite{DENISOVA} for a discussion of physical free boundary conditions). 
Starting from the nowadays classical local existence results \cite{SOLONNIKOV}, \cite{BEALE}, global existence was achieved in \cite{SYLVESTER}, \cite{TANAKA} for the scenario of an almost horizontal viscous fluid lying above a bottom and below vacuum. See also \cite{HATAYA}, \cite{GUO1}, \cite{GUO2}, \cite{GUO3} where the decay rate in time of the solution is studied in the previous situation in order to understand the long-time dynamics. In the two-fluid case (also known as internal waves problem in this scenario) global well-posedness and decay have been shown in \cite{WANGTICE}. See also \cite{MASMOUDI} for the vanishing viscosity limit problem for the free-boundary Navier-Stokes equations.

Recently several contributions have been made in the two-fluid case without viscosity jump with low regular positive density. First, Danchin and Mucha \cite{DANCHINMUCHA}, \cite{DANCHINPIECEWISE} showed the global well-posedness of \eqref{INH}-\eqref{NavierStokes} for initial densities allowing  discontinuities across $C^1$ interfaces
with a sufficiently small jump and small initial velocity in the Besov space $B_{p,1}^{2/p-1}$, $p\in[1,4)$ (see Section $2$ for the definition). For densities close enough to a positive constant and initial velocity in $\dot{B}_{p,1}^{2/p-1}\cap \dot{B}_{p,1}^{2/p-1+\eps}$, Huang, Paicu and Zhang \cite{HUANGPAICU} obtained solutions with $C^{1+\eps}$ flow  for small enough $\eps>0$.  Later Paicu, Zhang and Zhang \cite{PAICUZZ} obtained the global-wellposedness with initial data $u_0\in H^s$, $s\in(0,1)$ and initial positive density bounded from below and above removing the smallness conditions.

Based on these results and using paradifferential calculus and the techniques of striated regularity, Liao and Zhang have recently proved the persistence of $W^{k,p}$ regularity, $k\geq 3$, $p\in(2,4)$, for initial patches of the form
$$\rho_0=\rho_1 1_{D_0}+\rho_2 1_{D_0^c},$$
first assuming $\rho_1,\rho_2>0$ close to each other \cite{LIAOZHANG1}, then for any pair of positive constants \cite{LIAOZHANG2}.
By Sobolev embedding, this means that the boundary of the patch must be at least in $C^{2+\gamma}$ for some $\gamma>0$. Using the well-posedness result in \cite{DANCHINMUCHA}, Danchin and Zhang \cite{DANCHINZHANG} have recently obtained the propagation of $C^{1+\gamma}$ patches for small jump and small $u_0$ (also for large $u_0$ but only locally in time).

In this paper we consider the 2D density patch problem for Navier-Stokes without any smallness condition on the initial data and without any restriction on the density jump. 

First, we show that initial $C^{1+\gamma}$ density patches preserve their regularity globally in time for any $\rho_1,\rho_2>0$ and any $u_0\in H^{\gamma+s}$, $s\in(0,1-\gamma)$. We note that the cancellation in the tangential direction to the patch is not needed to propagate low regularities. Although one cannot expect to get the needed regularity for the velocity in Sobolev spaces, we will take advantage of the fact that $\rho$ remains as a patch with Lipschitz boundary. The quasilinear character of the coupling between density and velocity makes it harder to propagate the regularity of the velocity and hence that of the patch. This extra difficulty, compared to the same problem for Boussinesq system \cite{GANCEDOEDU}, is overcome by noticing that the characteristic function of a Lipschitz patch belongs to the multiplier space $\mathcal{M}(\dot{B}_{\infty,\infty}^{-1+\gamma})$ (see section below). Hence we will prove that $u\in L^1(0,T;C^{1+\gamma})$ and thus the propagation follows by the particle trajectories system
\begin{equation*}\label{trajectories}
\left\{\begin{aligned}
\frac{dX}{dt}(x,t)&=u(X(x,t),t),\\
X(x,0)&=x.
\end{aligned}\right.
\end{equation*}

Without considering regularity in the tangential direction to the density patch for the initial velocity, the initial conditions in \cite{LIAOZHANG2}, \cite{DANCHINZHANG} are  at the level of  $u_0\in B_{2,1}^{1+\epsilon}$ ($\epsilon>0$), $u_0\in B_{2,1}^{\gamma}$, respectively. Indeed, as in \cite{GANCEDOEDU}, from the results of maximum regularity of the linear heat equation, we deem $u_0\in H^{\gamma+s}$ is sharp at the scale of Sobolev spaces from this approach.

This low regularity result combined with new ideas allows us to show that the curvature of patches with initial $W^{2,\infty}$ regularity remains bounded for all time. 
Following the particle trajectory method to preserve the regularity, the curvature is controlled once that $\nabla^2u \in L^1(0,T; L^\infty)$. This is critical because at this level of derivatives the step function $\rho$ appears together with the nonlinearity. So in principle one could find that $\nabla^2u \in L^1(0,T; BMO)$. It is possible to use time weighted energy estimates, introduced in \cite{HOFF1}, \cite{HOFF2} for the compressible model and in \cite{PAICUZZ} for the incompressible case, combined with the characterization of a patch as a Sobolev multiplier to get higher regularity. In particular, to deal with the regularity of $u_t$ the convective derivative approach in \cite{LIAOZHANG2} can be used. Going further, the $C^{1+\gamma}$ regularity result in conjunction with the cancellation of singular integrals acting on low regular quadratic and cubic terms allows us to bootstrap to achieve the control of the evolution of the curvature ($\nabla^2 u \in L^1(0,T; L^\infty)$).

Finally, we continue the bootstrapping process to show a new proof for the propagation of regularity with initial $C^{2+\gamma}$ patches. Describing the dynamics of the patch by a level set, one can get advantage of the extra regularity in the tangential direction. In particular, in checking the evolution of this extra regularity one just needs to control the tangential direction of $\nabla^2 u$. Exploiting the smoothing properties of the Newtonian potential and the persistence of the regularity for the curvature, we are able to prove the propagation of $C^{2+\gamma}$ regularity. We realize that it is possible to find that extra cancellation dealing directly with singular integral operators.

\vspace{0.3cm}
	
The structure of the paper is as follows: In next section some definitions and notation that will be used along the paper are stated. Then, in Section 3 we show that the weak formulation we use to understand the solutions satisfies indeed the expected physical conditions at the interface. In Section 4 we prove the persistence of regularity for $C^{1+\gamma}$ patches and $u_0\in H^{\gamma+s}$. In Section 4 we show further that the curvature of the patches remains bounded. Finally, in Section 5 a proof of the propagation of $C^{2+\gamma}$ regularity in which we deal directly with the explicit expression of the tangential second derivatives of $u$ is given.

\section{Notation}

We include here some definitions and notations used along the paper.

We will denote by $\heatinv f$ the solution of the linear heat equation with force $f$ and zero initial condition:
\begin{equation*}\label{inverseheatequation}
(\partial_t-\Delta)^{-1}_0f:=\int_0^t e^{(t-\tau)\Delta}f(\tau)d\tau.
\end{equation*}
Above we use the standard notation $e^{t\Delta}f=\mathcal{F}^{-1}(e^{-t|\xi|^2}\hat{f})$, where $\hat{}$ and $\mathcal{F}^{-1}$ denote Fourier transform and its inverse.

The material derivate will be denote by $D_t f =\partial_tf + u\cdot \nabla f$.

We recall the definition of Besov spaces (see \cite{PDEBOOK} for details). Consider the following Littlewood-Paley decomposition: let $C=\{|\xi|\in\R^2:3/4\leq|\xi|\leq8/3\}$, fix a smooth radial function $\varphi$ supported in $C$ and satisfying 
\begin{equation*}\label{littlewoodp}
\begin{aligned}
\sum_{j\in\Z}\varphi(2^{-j}\xi)&=1,\hspace{0.3cm}\forall \xi\in \R^2\setminus\{0\}.
\end{aligned}
\end{equation*}
The homogeneous dyadic blocks are defined as $\Delta_jf=\mathcal{F}^{-1}(\varphi(2^{-j}\xi)\hat{f}(\xi))$ for $j\in \Z$.
Then, the homogeneous Besov spaces  $\dot{B}_{p,q}^\gamma(\R^2)$, $\gamma\in\R$, $p,q\in [1,\infty]$ are defined by
$$\dot{B}_{p,q}^\gamma(\R^2)=\{f\in S'(\R^2): \|f\|_{\dot{B}_{p,q}^\gamma}=\|2^{j\gamma}\|\Delta_j f\|_{L^p}\|_{l^q(\Z)}<\infty\},$$
where $S'(\R^2)$ denotes the space of tempered distributions over $\R^2$. We recall that $\dot{H}^s=\dot{B}_{2,2}^s$ and $\dot{C}^{k+\gamma}=\dot{B}_{\infty,\infty}^{k+\gamma}$ for $s\in \R$, $\gamma\in(0,1)$, $k\in \N\cup\{0\}$.

Let $E$ be a Banach space embedded in $S'(\R^n)$, $n=1,2$. We will use the spaces $L^p(0,T;E)$ with norm $\|f\|_{L^p_T(E)}:=\|\|f\|_{E}\|_{L^p(0,T)}$. 

We will say that $\partial D \in C([0,T]; E)$ if there exists a parametrization of the boundary 
\begin{equation}\label{parametrization}
\partial D(t) = \left\{z(\alpha,t)\in\R^2, \alpha\in[0,1]\right\}
\end{equation}
with $z \in C([0,T]; E)$.

 The {\em{multiplier space}} $\mathcal{M}(E)$ is the set of functions $\varphi$ such that $\varphi f\in E$ for all $f\in E$ and $$\|\varphi\|_{\mathcal{M}(E)}:=\sup_{\|f\|_{E}\leq 1}\|\varphi f\|_{E}<\infty.$$

\section{Weak solutions and physical conditions}

In this section we first state the definition of weak solutions for the system \eqref{INH}-\eqref{NavierStokes}. Later we show that under suitable regularity assumptions these solutions satisfy the naturally expected physical conditions (see e.g. \cite{DENISOVA} and the references therein):

\begin{itemize}
	\item[-] The interface moves with the fluid (no mass transfer):
	\begin{equation}\label{cond1}
	\begin{aligned}
	z_t(\alpha,t)\cdot n(\alpha,t)&=u\left(z(\alpha,t),t\right)\cdot n(\alpha,t),\\ 
	z(\alpha,0)&=z_0(\alpha),
	\end{aligned}
	\end{equation}
	where $z_0$ is a parametrization of the boundary of the patch $\partial D_0$.
	
	\item[-] Continuity of the velocity at the interface:
	\begin{equation}\label{cond2}
	 [u]|_{\partial D(t)}:= \lim_{\begin{subarray}
		{1}x\to x_0\in \partial D(t),\\\hspace{0.4cm}x\in D(t) 
		\end{subarray}} u(x)-\lim_{\begin{subarray}
		{1}x\to x_0\in \partial D(t),\\\hspace{0.15cm}x\in D(t)^c 
		\end{subarray}} u(x)=0,
	\end{equation}
	
	\item[-] Continuity of the stress tensor at the interface:
	\begin{equation}\label{cond3}
	\left[\mathcal{T}\cdot n\right]|_{\partial D(t)}=0,
	\end{equation}
	where $\mathcal{T}=-p\textup{ I}_\textup{d}+\left(\nabla u+\nabla u^*\right)$ and $\nabla u^*$ denotes transpose of $\nabla u$.
\end{itemize}

\begin{defi} We say that $(\rho, u, p)$ is a weak solution of the system \eqref{INH}-\eqref{NavierStokes} provided that
	\begin{equation}\label{weakincomp}
	\int_0^T\int_{\R^2} \nabla\varphi \cdot u \hspace{0.05cm}dx dt=0,\hspace{1cm}\forall \varphi \in C_c^\infty([0,T);C_0^\infty(\R^2)),
	\end{equation}
	and that for all $\phi$, $\Psi \in C^\infty_c([0,T);C_0^\infty(\R^2))$
	\begin{equation}\label{masscons}
	\int_{\R^2}\Psi(0)\rho_0dx+\int_0^T \int_{\R^2} \rho D_t\Psi\hspace{0.02cm} dx\hspace{0.02cm} dt=0,
	\end{equation}
	\begin{equation}\label{momentumcons}
	\!\int_{\R^2}\phi(0)\rho_0u_0\hspace{0.05cm}dx+\!\int_0^T\!\!\!\int_{\R^2}D_t \phi \cdot \rho u \hspace{0.05cm}dx\hspace{0.02cm} dt-\!\int_0^T\!\!\!\int_{\R^2}\nabla\phi \cdot\left(\nabla u+\nabla u^*\right)\hspace{0.05cm} dx \hspace{0.02cm}dt+\!\int_0^T\!\!\! \int_{\R^2} p\nabla\cdot \phi \hspace{0.05cm}dx\hspace{0.02cm}dt=0.
	\end{equation}
\end{defi}	
\vspace{0.5cm}
\begin{prop} Let $\left(\rho, u, p\right)$ be a weak solution of \eqref{INH}-\eqref{NavierStokes} with initial data as in Theorem \ref{Case1}. Then, the conditions \eqref{cond1}-\eqref{cond3} hold.	
\end{prop}
	
Proof:	The weak incompressibility condition \eqref{weakincomp} implies the continuity of the normal velocity at the interface, which jointly to the mass conservation \eqref{masscons} yields the interface dynamics condition \eqref{cond1} (see e.g. \cite{GANCEDO1}). Moreover, the results given in Theorem \ref{Case1} for initial data in $H^{\gamma+s}$ gives that $u\in L^1(0,T;C^{1+\gamma})$, hence the velocity is continuous also in the tangential direction for a.e. $t> 0$.  

We show then the continuity of the stress tensor. We can write
	\begin{equation}
	\begin{aligned}
	0&=\int_{\R^2}\phi(0)\rho_0u_0\hspace{0.05cm}dx+ \int_0^T\int_{D(t)\cup D(t)^c} D_t \phi\cdot \rho u\hspace{0.05cm}dxdt\\
	&\quad-\!\int_0^T\!\!\int_{D(t)\cup D(t)^c}\! \nabla\phi \cdot\left(\nabla u\!+\!\nabla u^*\right) dxdt + \int_0^T\int_{D(t)\cup D(t)^c} p \hspace{0.02cm}\nabla\cdot\phi dxdt.
	\end{aligned}
	\end{equation}
	Taking into account that the normal velocity is continuous at the interface and \eqref{cond1}, the regularity provided by Theorem \ref{Case1} allows us to integrate by parts to find that 
	\begin{equation}
	\begin{aligned}
	0&=\!\! \int_0^T\!\!\!\!\int_{D(t)}\!\phi\cdot \left(-\rho_1D_tu\!+\!\Delta u\!-\!\nabla p\right) \hspace{0.02cm}dxdt+\!\!\int_0^T\!\!\!\! \int_{D(t)^c} \!\phi \cdot \left(-\rho_2 D_tu\!+\!\Delta u\!-\!\nabla p\right)\hspace{0.02cm}dxdt\\	
	&\quad -\int_0^T\int_{\partial D(t)} \phi\hspace{0.05cm} n\cdot \left(\nabla u_1+\nabla u_1^*\right) d\sigma +\int_0^T\int_{\partial D(t)} \phi \hspace{0.05cm}n\cdot \left(\nabla u_2+\nabla u_2^*\right) d\sigma \\
	&\quad+\int_0^T\int_{\partial D(t)} \phi \hspace{0.05cm}p_1 n \hspace{0.05cm}d\sigma -\int_0^T\int_{\partial D(t)}\phi \hspace{0.05cm}p_2\hspace{0.02cm} n\hspace{0.05cm}d\sigma.
	\end{aligned}
	\end{equation}
	Thus we deduce that	
	\begin{equation}
	\begin{aligned}
	\int_0^T\int_{\partial D(t)}\phi \left[\left(p_1-p_2\right)\textup{I}_{\textup{d}}-\left( \left(\nabla u_1+\nabla u_1^*\right)-\left(\nabla u_2+\nabla u_2^*\right)\right)\right]\cdot n d\sigma =0.
	\end{aligned}
	\end{equation}
	and hence, as $p, u$ are regular enough \eqref{regularityp}, we conclude
	\begin{equation}
	\left[\left(-p\textup{ I}_{\textup{d}}+\left(\nabla u+\nabla u^*\right) \right)\cdot n\right]|_{\partial D(t)}=0,\hspace{0.5cm} \textup{a.e.} \hspace{0.2cm}t\in (0,T).
	\end{equation}

\qed

\section{Persistence of $C^{1+\gamma}$ regularity}\label{sec:2}

We present below the theorem that establishes the propagation of regularity for $C^{1+\gamma}$ patches in the case of positive density.

\vspace{0.3cm}
\begin{thm}\label{Case1}
Assume $\gamma\in(0,1)$, $s\in(0,1-\gamma)$, $\rho_1, \rho_2>0$. Let $D_0\subset \R^2$ be a bounded simply connected domain with boundary $\partial D_0\in C^{1+\gamma}$, $u_0\in H^{\gamma+s}$ a divergence-free vector field,  
\begin{equation*}
\rho_0(x)=\rho_1 1_{D_0}(x)+\rho_2 1_{D_0^c}(x),
\end{equation*}
and $1_{D_0}$ the characteristic function of $D_0$. Then, there exists a unique global solution $(u,\rho)$ of \eqref{INH}-\eqref{NavierStokes} such that
$$\rho(x,t)=\rho_1 1_{D(t)}(x)+\rho_2 1_{D(t)^c}(x) \hspace{0.2cm}{\rm{and}} \hspace{0.2cm} \partial D\in C([0,T]; C^{1+\gamma}),$$
where $D(t)=X(D_0,t)$ with $X$ the particle trajectories associated to the velocity field.\newline
Moreover, 
\begin{equation*}
\begin{aligned}
u&\in C([0,T]; H^{\gamma+s})\cap L^1(0,T;C^{1+\gamma+\tilde{s}}),\\
t^{\frac{1-(\gamma+s)}{2}}u &\in L^\infty(0,T;H^1),\hspace{0.1cm} t^{\frac{2-(\gamma+s)}{2}}u\in L^\infty(0,T;H^2),\\
t^{\frac{2-(\gamma+s)}{2}} u_t&\in L^\infty([0,T];L^2)\cap L^2([0,T];H^1),
\end{aligned}
\end{equation*}
for any $T>0$, $\tilde{s}\in (0,s)$.  
\end{thm}

Proof: First, as $u_0\in H^{\gamma+s}$ and $0<\min\{\rho_1,\rho_2\}<\rho_0<\max\{\rho_1,\rho_2\}<\infty$, the results in \cite{PAICUZZ} yield the following estimates for any $T\geq 0$:
\begin{equation}\label{apriori}
\begin{aligned}
A_0(T)&\leq C(\|u_0\|_{L^2}),\\
 A_1(T)&\leq C(\|u_0\|_{H^{\gamma+s}}),\\
A_2(T)&\leq C(\|u_0\|_{H^{\gamma+s}}),\\
\int_0^T \|\nabla u \|_{L^\infty}dt&\leq C(T,\|u_0\|_{H^{\gamma+s}}),
\end{aligned}
\end{equation}
where the constant $C$ also depends on $\rho_1, \rho_2$, and $A_0, A_1, A_2$ are defined by
\begin{equation*}
\begin{aligned}
A_0(T)&=\sup_{[0,T]}\|\sqrt{\rho}u\|_{L^2}^2+\int_0^T\|\grad u\|_{L^2}^2dt,\\
A_1(T)&=\sup_{[0,T]}t^{1-(\gamma+s)}\|\nabla u\|_{L^2}^2,\\
A_2(T)&=\sup_{[0,T]}t^{2-(\gamma+s)}\left( \|\sqrt{\rho}u_t\|_{L^2}^2+\|\Delta u\|_{L^2}^2+\|\nabla p\|_{L^2}^2 \right)+\int_0^T t^{2-(\gamma+s)}\|\grad u_t\|_{L^2}^2.
\end{aligned}
\end{equation*}
In particular, we note that by interpolation we get
\begin{equation*}
\|u\|_{H^{1+\gamma+\tilde{s}}}\leq \|u\|_{H^1}^{1-\gamma-\tilde{s}}\|u\|_{H^2}^{\gamma+\tilde{s}}\leq c t^{-\frac{1-\gamma-s}{2}(1-\gamma-\tilde{s})-\frac{2-\gamma-s}{2}(\gamma+\tilde{s})},
\end{equation*}
and therefore
\begin{equation}\label{uL2}
u\in L^p(0,T;H^{1+\gamma+\tilde{s}}),\hspace{1cm} p\in[1,2/(1-(s-\tilde{s}))).
\end{equation}
Proceeding by interpolation again, it follows that
\begin{equation*}
\int_0^T \!\!\|u_t\|_{H^{\gamma+\tilde{s}}}^q dt\leq\!\!\! \int_0^T\!\! \left(\|u_t\|_{L^2}^{1-\gamma-\tilde{s}}\|u_t\|_{H^1}^{\gamma+\tilde{s}}\right)^qdt\!\leq \!c\!\! \int_0^T\!\!\!t^{-\frac{2-\gamma-s}{2}(1-\gamma-\tilde{s})q}\frac{\|u_t\|_{H^1}^{q(\gamma+\tilde{s})}t^{\frac{2-\gamma-s}{2}(\gamma+\tilde{s})q}}{t^{\frac{2-\gamma-s}{2}(\gamma+\tilde{s})q}}dt,
\end{equation*}
hence by H\"older inequality we conclude 
\begin{equation}\label{utHs}
u_t\in L^q(0,T;H^{\gamma+\tilde{s}}), \hspace{1cm}q\in[1,2/(2-(s-\tilde{s}))).
\end{equation}

Next, we rewrite \eqref{NavierStokes} as a forced heat equation
\begin{equation*}
u_t-\Delta u=-\rho u\cdot\nabla u+(1-\rho)u_t-\nabla p.
\end{equation*}
We apply first the Leray projector $\mathbb{P}=\textup{I}_\textup{d}-\nabla\Delta^{-1}(\nabla\cdot \hspace{0.1cm})$ to obtain
\begin{equation}\label{forceheatleray}
u_t-\Delta u=-\mathbb{P} (\rho  u\cdot\nabla u))+\mathbb{P}((1-\rho)u_t),
\end{equation}
and denote
\begin{equation}\label{decomposition}
\begin{gathered}
u=v_1+v_2+v_3,\\
v_1=e^{t\Delta}u_0,\hspace{1cm} v_2=-\heatinv \mathbb{P}(\rho u\cdot\nabla u),\hspace{1cm} v_3=\heatinv \mathbb{P}((1-\rho)u_t).
\end{gathered}
\end{equation}
Recalling the following particular case of Gagliardo-Nirenberg inequality
\begin{equation}
\label{gagliardonirenberg}
\|f\|_{L^r(\R^2)}\leq c\|f\|_{L^2}^{2/r} \|\nabla f\|_{L^2}^{1-2/r},\hspace{0.7cm}r\in[2,\infty),
\end{equation}
we deduce that
\begin{equation*}
\|u\cdot \nabla u\|_{L^2}\leq \|u\|_{L^2}^{1/2}\|\nabla u\|_{L^2}\|\nabla^2u\|_{L^2}^{1/2}\leq c t^{-1+\frac34(\gamma+s)}.
\end{equation*}
Then, the splitting in \eqref{decomposition} provides 
\begin{equation*}
\|u\|_{L^\infty_T(H^{\gamma+s})}\leq c\|u_0\|_{H^{\gamma+s}}+\|\left(1-\Delta\right)^{\frac{\gamma+s}{2}}(v_2+v_3)\|_{L^\infty(L^2)},
\end{equation*} 
so that using the decay properties of the heat kernel (see e.g. Appendix D in \cite{RODRIGO}) and Young's inequality for convolutions we obtain
\begin{equation}
\begin{aligned}
\|u\|_{L^\infty_T(H^{\gamma+s})}&\leq c\|u_0\|_{H^{\gamma+s}}+c\left(\|u_0\|_{H^{\gamma+s}}\right)\left|\left| \int_0^t  \|\left(1-\Delta\right)^{\frac{\gamma+s}{2}}K(t-\tau)\|_{L^1} \tau^{-1+\frac{\gamma+s}{2}}  d\tau\right|\right|_{L^\infty_T}\\
&\leq c\|u_0\|_{H^{\gamma+s}}+c\left(\|u_0\|_{H^{\gamma+s}}\right) \left|\left| \int_0^t (t-\tau)^{-\frac{\gamma+s}{2}} \tau^{-1+\frac{\gamma+s}{2}}d\tau\right|\right|_{L^\infty_T}\leq c\left(\|u_0\|_{H^{\gamma+s}}\right).
\end{aligned}
\end{equation}

Notice that from \eqref{apriori} the velocity field satisfies $u\in L^1(0,T; W^{1,\infty})$, so the initial density is transported and remains as a patch
\begin{equation*}
\rho(x,t)=\rho_1 1_{D(t)}(x)+\rho_2 1_{D(t)^c}(x)
\end{equation*}
with Lipschitz boundary. It is known (see \cite{TRIEBEL}) that the characteristic function of a Lipschitz bounded domain belongs to the multiplier space $\mathcal{M}(\dot{B}_{a,b}^s)$ if and only if $-1+\frac{1}{a}<s<\frac{1}{a}$, where $a,b\in [1,\infty]$. Therefore we have that
\begin{equation}\label{densitymultipliergeneral}
\rho \in L^\infty(0,T;\mathcal{M}(\dot{B}_{a,b}^{s})),\hspace{0.5cm} -1+\frac{1}{a}<s<\frac{1}{a},
\end{equation}
so in particular
\begin{equation}\label{densitymultiplier}
\rho \in L^\infty(0,T;\mathcal{M}(\dot{B}_{\infty,\infty}^{-1+\gamma+s})).
\end{equation}
We write \eqref{NavierStokes} as follows
\begin{equation}\label{laplace}
\Delta u= \mathbb{P}\left(\rho D_t u\right).
\end{equation}
From \eqref{uL2} we deduce  $\nabla\cdot(u\otimes u)\in L^{p/2}(0,T;H^{\gamma+\tilde{s}})$. This joined to \eqref{utHs} yields that $D_t u\in L^q(0,T; H^{\gamma+\tilde{s}})\hookrightarrow L^q(0,T; \dot{H}^{\gamma+\tilde{s}})$ for $q\in[1,2/(2-(s-\tilde{s})))$. By embedding in Besov spaces we have that
$$D_t u\in L^q(0,T; \dot{B}_{\infty,\infty}^{-1+\gamma+\tilde{s}}),$$
and taking into account \eqref{densitymultiplier} it follows that
$$\rho D_t u\in L^q(0,T; \dot{B}_{\infty,\infty}^{-1+\gamma+\tilde{s}}).$$
Finally, recalling that the Leray projector is a Fourier multiplier of degree zero and hence it is bounded in H\"older spaces, we take the inverse of the Laplacian in \eqref{laplace} to find that
\begin{equation*}
\|u\|_{\dot{C}^{1+\gamma+\tilde{s}}}\leq c\|\Delta^{-1}\rho D_t u\|_{\dot{C}^{1+\gamma+\tilde{s}}}\leq c \|\rho D_t u\|_{\dot{B}_{\infty,\infty}^{-1+\gamma+\tilde{s}}}.
\end{equation*}
From \eqref{uL2} it is clear that $u\in L^q(0,T;L^\infty)$, therefore we conclude that $u\in L^q(0,T; C^{1+\gamma+\tilde{s}})$. Hence, using the particle trajectories system \eqref{trajectories} and Gr\"onwall's inequality we obtain
\begin{equation*}
\|\grad X\|_{C^\gamma}\leq  \|\grad X_0\|_{C^\gamma}e^{\int_0^t\|\grad u\|_{L^\infty}d\tau}+\int_0^t\|\grad u(\tau)\|_{C^\gamma}\|\grad X(\tau)\|_{L^\infty}^{1+\gamma}e^{\int_\tau^t\|\grad u\|_{L^\infty}ds},
\end{equation*}
which yields the persistence of $C^{1+\gamma}$ regularity of the density patch $\|z\|_{L^\infty(0,T;C^{1+\gamma})}\leq C(T).$

\begin{rem} From the momentum equation it is easy to see that
	\begin{equation*}
	\|\Delta u\|_{H^\mu}^2+\|\nabla p\|_{H^\mu}^2= \|\rho D_tu \|_{H^\mu}^2.
	\end{equation*}	
Noticing that $\rho \in \mathcal{M}(H^{\sigma})$, $\sigma\in (-1/2,1/2)$ and recalling that $D_t u\in L^q(0,T;H^{\gamma+s})$, it follows that
\begin{equation}\label{regularityp}
 p\in L^q(0,T;\dot{H}^{1+\mu}),\hspace{0.5cm} u\in L^q(0,T;H^{2+\mu}),\hspace{0.5cm}\mu<\min\{1/2,\gamma+s\}.
\end{equation}
\end{rem}

\qed

\section{Persistence of $W^{2,\infty}$ regularity}\label{sec:3}

In this section we show that the curvature of the patch is bounded for all time if initially has $W^{2,\infty}$ boundary.

\begin{thm}\label{Case2}
	Assume $s\in(0,1)$, $\rho_1, \rho_2>0$. Let $D_0\subset \R^2$ be a bounded simply connected domain with boundary $\partial D_0\in W^{2,\infty}$, $u_0\in H^{1+s}$ a divergence-free vector field,  
	\begin{equation*}
	\rho_0(x)=\rho_1 1_{D_0}(x)+\rho_2 1_{D_0^c}(x),
	\end{equation*}
	and $1_{D_0}$ the characteristic function of $D_0$. Then, there exists a unique global solution $(u,\rho)$ of \eqref{INH}-\eqref{NavierStokes} such that
	$$\rho(x,t)=\rho_1 1_{D(t)}(x)+\rho_2 1_{D(t)^c}(x) \hspace{0.2cm}{\rm{and}} \hspace{0.2cm} \partial D\in C([0,T];W^{2,\infty}),$$
	where $D(t)=X(D_0,t)$ with $X$ the particle trajectories associated to the velocity field.\newline
	Moreover, 
	\begin{equation*}
	\begin{aligned}
	u&\in C([0,T];H^{1+s})\cap L^2(0,T;H^{2+\mu})\cap L^p(0,T;W^{2,\infty}),\\
	t^{\frac{1-s}{2}} u_t&\in L^\infty(\left[0,T\right];L^2)\cap L^2(\left[0,T\right];H^1),\\
	t^{\frac{2-s}{2}} D_tu&\in L^\infty(\left[0,T\right];H^1)\cap L^2(\left[0,T\right];H^2),
	\end{aligned}
	\end{equation*}
	for any $T>0$, $\mu<\min\{1/2,s\}$ and $p\in[1,2/(2-s))$. If $s<1/2$ it also holds that $u\in L^q(0,T;H^{2+\delta})$ for any $\delta\in(s,1/2)$ with $q\in[1,2/(1+\delta-s))$.

\end{thm}
	
Proof: First, we notice that once we get $u\in L^p(0,T; W^{2,\infty})$ the propagation of regularity for the patch follows by considering the particle trajectories associated to the flow. Hence, we proceed to prove that the velocity belongs to that space.

\vspace{0.1cm}
\subsection{Regularity of $u_t$}
\vspace{0.3cm}

We start by proving in this section that $\displaystyle t^{\frac{1-s}{2}} u_t\in L^\infty(\left[0,T\right];L^2)\cap L^2(\left[0,T\right];H^1)$.
As before, it is easy to get
\begin{equation}\label{balanceL2}
\|\sqrt{\rho} u\|_{L^2}^2+\int_0^t \|\nabla u\|_{L^2}^2d\tau\leq  c\hspace{0.02cm}\|u_0\|_{L^2}^2.
\end{equation}
We now take inner product with $u_t$ and use Young's inequality to obtain
\begin{equation}\label{balance1}
\|\sqrt{\rho}u_t\|_{L^2}^2+\frac{d}{dt}\|\nabla u\|_{L^2}^2\leq c\hspace{0.02cm}\|u\|_{L^4}^2\|\nabla u\|_{L^4}^2.
\end{equation}
Using \eqref{gagliardonirenberg} with $r=4$, from the velocity equation and Young's inequality one infers that
\begin{equation}\label{Stokespri}
\|\nabla^2 u\|_{L^2}^2+\|\nabla p\|_{L^2}^2\leq c \left(\|\sqrt{\rho}u_t\|_{L^2}^2+\|u\|_{L^2}^2\|\nabla u\|_{L^2}^4\right).
\end{equation}
Hence, applying \eqref{gagliardonirenberg} again in \eqref{balance1} and using \eqref{Stokespri}, we get
\begin{equation*}
\|\sqrt{\rho}u_t\|_{L^2}^2+\frac{d}{dt}\|\nabla u\|_{L^2}^2\leq c\hspace{0.02cm}\|u\|_{L^2}^2\|\nabla u\|_{L^2}^4,
\end{equation*}
so by Gr\"onwall's inequality and \eqref{balanceL2} we conclude that
\begin{equation*}
\|\nabla u\|_{L^2}^2+\int_0^t \|\sqrt{\rho}u_t\|_{L^2}^2d\tau\leq \|\nabla u_0\|_{L^2}^2e^{c\hspace{0.02cm}\|u_0\|_{L^2}t}.
\end{equation*}
We can close the estimates for $u$ and $u_t$ at this level of regularity:
\begin{equation}\label{lionsestimates}
\begin{aligned}
\|u\|_{L^\infty_T(L^2)}+\|u\|_{L^2_T(H^1)}&\leq c\hspace{0.02cm}\|u_0\|_{L^2},\\
\|u\|_{L^\infty_T(H^1)}+\|u\|_{L^2_T(H^2)}+\|u_t\|_{L^2_T(L^2)}&\leq c\left(\|u_0\|_{L^2}\right)\|u_0\|_{H^1}.
\end{aligned}
\end{equation}
From this last estimate, \eqref{Stokespri} rewrites as
\begin{equation}\label{Stokes}
\|\nabla^2 u\|_{L^2}^2+\|\nabla p\|_{L^2}^2\leq c\hspace{0.02cm} \left(\|\sqrt{\rho}u_t\|_{L^2}^2+1\right).
\end{equation}
We proceed next by an interpolation argument. First, we consider the linear momentum equation for $v$
\begin{equation}\label{vequation}
\rho v_t+\rho u\cdot\nabla v-\Delta v+\nabla p=0,\hspace{1cm} \rho_t+u\cdot \nabla\rho=0.
\end{equation}
By previous arguments it follows that 
\begin{equation}\label{vL2balance}
\|v\|_{L^\infty_T(H^1)}^2+\|\sqrt{\rho}v_t\|_{L^2_T(L^2)}^2\leq c\left(\|u_0\|_{L^2}\right)\|v_0\|_{H^1}^2,
\end{equation}
\begin{equation*}
\|\nabla^2 v\|_{L^2}^2\leq c\hspace{0.02cm}\left( \|\sqrt{\rho}v_t\|_{L^2}^2+\|u\|_{L^2}^2\|\nabla u\|_{L^2}^2\|\nabla v\|_{L^2}^2 \right),
\end{equation*}
and hence
\begin{equation}\label{vStokes}
\|\nabla^2 v\|_{L^2}^2\leq c\hspace{0.02cm}\left( \|\sqrt{\rho}v_t\|_{L^2}^2+\|v_0\|_{H^1}^2 \right),
\end{equation}
Derivation in time of \eqref{vequation} yields the following equation
\begin{equation*}
\rho v_{tt}+\rho u\cdot\nabla v_t-\Delta v_t+\nabla p_t=-\rho_t v_t-\rho_t u\cdot\nabla v-\rho u_t\cdot\nabla v,
\end{equation*}
and thus we obtain 
\begin{equation*}
\frac12\frac{d}{dt}\|\sqrt{\rho}v_t\|_{L^2}^2+\|\nabla v_t\|_{L^2}^2=-\int_{\R^2}\rho_t |v_t|^2dx-\int_{\R^2}\rho_t u\cdot\nabla v\cdot v_t dx-\int_{\R^2}\rho u_t\cdot\nabla v\cdot v_t dx,
\end{equation*}
where we have used that $\rho_t=-u\cdot\nabla\rho$. Multiplication by the weight $t$ and integration in time implies that
\begin{equation}\label{integratedenergybalancedt}
\frac{t}{2}\|\sqrt{\rho}v_t\|_{L^2}^2+\int_0^t\tau \|\nabla v_t\|_{L^2}^2 d\tau=I_1+I_2+I_3+I_4,
\end{equation}
where
\begin{equation*}
\begin{aligned}
I_1&=\frac12 \int_0^t \|\sqrt{\rho}v_t\|_{L^2}^2d\tau, \hspace{2.6cm}I_2=-\int_0^t\tau\int_{\R^2} \rho_t|v_t|^2dx d\tau,  \\
I_3&=-\int_0^t \tau\int_{\R^2} \rho_t u\cdot\nabla v\cdot v_t dx d\tau,\hspace{1cm}I_4=-\int_0^t \tau\int_{\R^2} \rho u_t\cdot\nabla v\cdot v_t dx d\tau.
\end{aligned}
\end{equation*}
The first term is controlled by \eqref{vL2balance} as follows
\begin{equation}\label{I1}
I_1\leq c\left(\|u_0\|_{L^2} \right)\|v_0\|_{H^1}^2.
\end{equation}
Recalling that $u$ is divergence-free and that $\rho_t=-u\cdot\nabla \rho$, integration by parts in $I_2$ yields the following
\begin{equation*}
I_2\leq \int_0^t \tau\int_{\R^2}\rho u\cdot\nabla |v_t|^2dx d\tau\leq c\hspace{0.02cm}\int_0^t\tau \|\nabla v_t\|_{L^2}\|u\|_{L^4}\|v_t\|_{L^4}d\tau.
\end{equation*}
By \eqref{gagliardonirenberg} and Young's inequality $I_2$ is bounded by
\begin{equation}\label{I2}
I_2\leq \frac{1}{10}\int_0^t \tau\|\nabla v_t\|_{L^2}^2d\tau+c(\|u_0\|_{L^2})\int_0^t \tau \|\sqrt{\rho}v_t\|_{L^2}^2\|\nabla u\|_{L^2}^2 d\tau.
\end{equation}
Using again \eqref{gagliardonirenberg} and \eqref{vStokes} we get that 
\begin{equation*}
\begin{aligned}
I_4&\leq c\left(\|u_0\|_{H^1}\right) \int_0^t \tau \|\sqrt{\rho}u_t\|_{L^2}\|\nabla v\|_{L^2}^{1/2}\|\sqrt{\rho}v_t\|_{L^2}\|\nabla v_t\|_{L^2}^{1/2}d\tau\\
&\quad+c\left(\|u_0\|_{H^1}\right)\int_0^t \tau\|\sqrt{\rho}u_t\|_{L^2}\|\nabla v\|_{L^2}\|\sqrt{\rho}v_t\|_{L^2}^{1/2}\|\nabla v_t\|_{L^2}^{1/2}d\tau,
\end{aligned}
\end{equation*}
and therefore by \eqref{vL2balance} and Young's inequality it follows that
\begin{equation}\label{I4}
\begin{aligned}
I_4&\leq \frac{1}{10}\int_0^t \tau \|\nabla v_t\|_{L^2}^2d\tau+c\left(\|u_0\|_{H^1}\right)\int_0^t\tau \|\sqrt{\rho}u_t\|_{L^2}^2\|\sqrt{\rho}v_t\|_{L^2}^2 d\tau\\
&\quad+c\left(\|u_0\|_{L^2},T\right)\|v_0\|_{H^1}^2.
\end{aligned}
\end{equation}
After integration by parts, $I_3$ is decomposed as follows
\begin{equation*}
I_3=I_{31}+I_{32}+I_{33},
\end{equation*}
where
\begin{equation*}
\begin{aligned}
I_{31}&=-\int_0^t \tau\int_{\R^2} \rho v_t\cdot\nabla u\cdot\nabla v\cdot u \hspace{0.05cm}dxd\tau,\\
I_{32}&=-\int_0^t \tau\int_{\R^2} \rho (u\otimes u):\nabla^2 v\cdot v_t d\tau dx,\\
I_{33}&=-\int_0^t \tau\int_{\R^2} \rho (u\cdot\nabla v)\cdot(u\cdot\nabla v_t) dxd\tau.
\end{aligned}
\end{equation*}
First we use again \eqref{gagliardonirenberg} and Young's inequality  to obtain
\begin{equation*}
\begin{aligned}
I_{31}&\leq \int_0^t\tau \|\rho u\cdot \nabla u\|_{L^2}\|v_t\|_{L^4}\|\nabla v\|_{L^4}d\tau\\
& \leq \frac{1}{10}\int_0^t \tau \|\nabla v_t\|_{L^2}^2d\tau+c\hspace{0.02cm}\int_0^t \tau\|u\cdot \nabla u\|_{L^2}^{4/3}\|\sqrt{\rho} v_t\|_{L^2}^{2/3}\|\nabla v\|_{L^2}^{2/3}\|\nabla^2 v\|_{L^2}^{2/3}d\tau.
\end{aligned}
\end{equation*}
After using \eqref{vStokes}, Young's inequality and the previous estimates on $u$ and $v$ yield
\begin{equation*}
\begin{aligned}
I_{31}&\leq \frac{1}{10}\int_0^t \tau\|\nabla v_t\|_{L^2}^2d\tau +c\left(\|u_0\|_{H^1}\right)\int_0^t \tau\|\sqrt{\rho}v_t\|_{L^2}^2\left(1+\|\nabla^2 u\|_{L^2}^2\right)d\tau\\
&\quad+c\left(\|u_0\|_{L^2},T\right)\|v_0\|_{H^1}^2.
\end{aligned}
\end{equation*}
Since \eqref{gagliardonirenberg} and \eqref{lionsestimates} give $\|u\|_{L^8}\leq c\left(\|u_0\|_{L^2}\right)\|u_0\|_{H^1}$, using \eqref{gagliardonirenberg}, \eqref{vStokes} and Young's inequality the terms $I_{32}$ and $I_{33}$ are bounded as $I_{31}$. Therefore, 
\begin{equation}\label{I3}
\begin{aligned}
I_3&\leq \frac{1}{10}\int_0^t \tau\|\nabla v_t\|_{L^2}^2d\tau +c\hspace{0.02cm}\int_0^t \tau\|\sqrt{\rho}v_t\|_{L^2}^2\left(1+\|\nabla^2 u\|_{L^2}^2\right)d\tau\\
&\quad+c\left(\|u_0\|_{L^2},T\right)\|v_0\|_{H^1}^2.
\end{aligned}
\end{equation}
Joining the above bounds \eqref{I1}-\eqref{I3} we get from \eqref{integratedenergybalancedt} that
\begin{equation*}
\begin{aligned}
t\|\sqrt{\rho}v_t\|_{L^2}^2+\int_0^t \tau \|\nabla v_t\|_{L^2}^2d\tau &\leq c\left( \|u_0\|_{L^2},T\right)\|v_0\|_{H^1}^2\\
&\quad+c\hspace{0.02cm}\int_0^t \tau \|\sqrt{\rho}v_t\|_{L^2}^2\left(\|\nabla^2 u\|_{L^2}^2+\|\sqrt{\rho}u_t\|_{L^2}^2+1\right)d\tau,
\end{aligned}
\end{equation*}
thus by Gr\"onwall's inequality we finally find
\begin{equation*}
t\|\sqrt{\rho}v_t\|_{L^2}^2+\int_0^t \tau \|\nabla v_t\|_{L^2}^2d\tau\leq c\left( \|u_0\|_{H^1},T\right)\|v_0\|_{H^1}^2.
\end{equation*}
We notice that from \eqref{vequation} we have $\|\sqrt{\rho} v_t\|_{L^2}(\tau)\leq c\left( \|u_0\|_{H^1}\right)\|v\|_{H^2}(\tau)$ for all $\tau\geq 0$, so if we assumed $v_0\in H^2$, repeating the steps above without weights would lead to
\begin{equation*}
\|\sqrt{\rho}v_t\|_{L^2}^2+\int_0^t \|\nabla v_t\|_{L^2}^2d\tau\leq c\left( \|u_0\|_{H^1},T\right)\|v_0\|_{H^2}^2.
\end{equation*}
Finally, by linear interpolation between the last two inequalities \cite{LUNARDI} we conclude that
\begin{equation*}
t^{1-s}\|\sqrt{\rho}u_t\|_{L^2}^2+\int_0^t \tau^{1-s}\|\nabla u_t\|_{L^2}^2d\tau \leq c\left( \|u_0\|_{H^1},T\right)\|u_0\|_{H^{1+s}}^2.
\end{equation*} 
Using \eqref{Stokes} we are able to finally find
\begin{equation}\label{apriori1}
t^{1-s}\|\sqrt{\rho}u_t\|_{L^2}^2+ t^{1-s}\|\nabla^2 u\|_{L^2}^2+\int_0^t \tau^{1-s}\|\nabla u_t\|_{L^2}^2d\tau \leq c\left( \|u_0\|_{H^1},T\right)\|u_0\|_{H^{1+s}}^2.
\end{equation}
We note that by Sobolev interpolation we have in particular that for $\tilde{s}\in(0,s)$
\begin{equation}\label{uths}
u_t\in L^2(0,T; H^{\tilde{s}}).
\end{equation}

\vspace{0.1cm}
\subsection{Higher regularity for $u$}
\vspace{0.3cm}

This section is devoted to prove that $\displaystyle u\in L^\infty(0,T; H^{1+s})\cap L^2(0,T; H^{2+\mu})$.
First, we notice that by Theorem \ref{Case1} $\rho(t)$ remains as a patch with Lipschitz boundary for all $t\geq 0$ and hence it is known that 
\begin{equation*}
\rho \in L^\infty(0,T;\mathcal{M}(H^\sigma))\hspace{0.5cm}\sigma\in\left(-\frac12,\frac12\right).
\end{equation*}
From this and the estimates \eqref{lionsestimates}, \eqref{uths} we infer that 
\begin{equation*}
(1-\rho) u_t, \hspace{0.2cm}\rho u\cdot\nabla u \in L^2(0,T;H^{\mu}),\hspace{0.5cm}\mu<\min\left\{\frac12,s\right\}.
\end{equation*}
Thus, by classical properties of the heat equation applied to \eqref{forceheatleray} we conclude that $u\in L^2(0,T;H^{2+\mu})$ with
\begin{equation}\label{uH2mu}
\|u\|_{L^2_T(H^{2+\mu})}\leq c\left(\|u_0\|_{H^{1+s}},T\right).
\end{equation}
To get $u\in L^\infty(0,T;H^{1+s})$ we will use that \eqref{apriori1} and \eqref{Stokes} implies 
\begin{equation}\label{aux}
t^{1-s}\|\sqrt{\rho} u_t\|_{L^2}^2+ t^{1-s}\|\sqrt{\rho} u\cdot\nabla u\|_{L^2}^2\leq  c\left( \|u_0\|_{H^{1+s}},T\right).
\end{equation}
Then, we write the solution of \eqref{forceheatleray} as
\begin{equation*}
u=e^{t\Delta}u_0+\heatinv \left(\mathbb{P}\left((1-\rho)u_t\right) -\mathbb{P}\left(\rho u\cdot\nabla u\right)\right),
\end{equation*}
thus we have that
\begin{equation*}
\begin{aligned}
\|u\|_{L^\infty_T(H^{1+s})}&\leq c\|u_0\|_{H^{1+s}}+c\left|\left| \left(1-\Delta\right)^{\frac{1+s}{2}}\heatinv \left(\mathbb{P}\left((1-\rho)u_t\right) -\mathbb{P}\left(\rho u\cdot\nabla u\right)\right)\right|\right|_{L^\infty_T(L^2)}.
\end{aligned}
\end{equation*}
Applying Young's inequality for convolution, \eqref{aux} and the decay properties of the heat kernel we get 
\begin{equation*}
\begin{aligned}
\|u\|_{L^\infty_T(H^{1+s})}&\leq c\|u_0\|_{H^{1+s}}+c\left|\left| \int_0^t\|\left(1-\Delta\right)^{\frac{1+s}{2}} K(t-\tau)\|_{L^1}\tau^{-\frac{1-s}{2}}  d\tau\right|\right|_{L^\infty_T}\\
&\leq c\|u_0\|_{H^{1+s}}+c\left|\left|\int_0^t (t-\tau)^{-\frac{1+s}{2}}\tau^{-\frac{1-s}{2}}  d\tau\right|\right|_{L^\infty_T},
\end{aligned}
\end{equation*}
so that we conclude
\begin{equation}\label{uLinf}
\|u\|_{L^\infty_T(H^{1+s})}\leq c\left( \|u_0\|_{H^{1+s}},T\right).
\end{equation}

\vspace{0.3cm}
\subsection{Higher regularity for $D_t u$}
\vspace{0.3cm}

We will show that $\displaystyle t^{\frac{2-s}{2}} D_t u\in L^\infty(\left[0,T\right];H^1)\cap L^2(\left[0,T\right];H^2)$.
Applying $D_t$ to \eqref{NavierStokes} yields 
\begin{equation}\label{Dt2}
\rho D_t^2 u-\Delta D_t u +\nabla D_t p=-2\nabla u_i\cdot \partial_i\nabla u+\Delta u\cdot\nabla u-\nabla u^T\cdot\nabla p,
\end{equation} 
where Einstein summation convention is used. By definition of $D_t u$, it follows directly from previous estimates \eqref{apriori1} and \eqref{uLinf} that
\begin{equation}\label{Dtuestimates}
t^{1-s}\|\sqrt{\rho}D_t u\|_{L^2}^2+\int_0^t \tau^{1-s}\|\nabla D_t u\|_{L^2}^2\leq c\left( \|u_0\|_{H^1},T\right)\|u_0\|_{H^{1+s}}^2.
\end{equation}
 In what follows we will denote 
 \begin{equation*}
 F=F(\nabla u, \nabla^2u, \nabla p)=-2\nabla u_i\cdot \partial_i\nabla u+\Delta u\cdot\nabla u-\nabla u^T\cdot\nabla p.
 \end{equation*}
 As $\|\nabla^2 u\|_{L^2}+\|\nabla p\|_{L^2}\leq \|\rho D_tu\|_{L^2}$, using \eqref{Dtuestimates} we notice that
 \begin{equation*}
 \|F\|_{L^2}^2\leq c\|\nabla u\|_{L^\infty}^2\|\rho D_t u\|_{L^2}^2\leq c\hspace{0.05cm}t^{-1+s}\|u\|_{H^{2+\eps}}^2,
 \end{equation*}
 so from \eqref{uH2mu} we find
 \begin{equation}\label{boundF} 
 \int_0^t \tau^{1-s}\|F\|_{L^2}^2 d\tau\leq c\left(\|u_0\|_{H^{1+s}},T\right).
\end{equation}
By taking dot product of \eqref{Dt2} with $D_t^2u$ and integrating in time we find
\begin{equation}\label{Dt2balance}
\frac12t^{2-s}\|\nabla D_t u\|_{L^2}^2+\int_0^t \tau^{2-s}\|\sqrt{\rho}D_t^2 u\|_{L^2}^2d\tau=L_1+L_2+L_3+L_4,
\end{equation}
where
\begin{equation*}
\begin{aligned}
L_1&=\frac{2-s}{2}\int_0^t \tau^{1-s}\|\nabla D_t u\|_{L^2}^2d\tau,\hspace{1.1cm}L_2=-\int_0^t \tau^{2-s}\int_{\R^2}D^2_t u\cdot\nabla D_t p \hspace{0.2cm}dxd\tau,\\
L_3&=\int_0^t \tau^{2-s}\int_{\R^2} D^2_t u\cdot F\hspace{0.05cm}dxd\tau,\hspace{1.2cm}L_4=\int_0^t \tau^{2-s}\int_{\R^2}\nabla D_t u: \nabla(u\cdot \nabla D_t u)dxd\tau.
\end{aligned}
\end{equation*}
The first term, $L_1$, is bounded by \eqref{Dtuestimates},
\begin{equation}\label{L1}
L_1\leq  c\left( \|u_0\|_{H^{1+s}}^2,T\right).
\end{equation}
 while using \eqref{boundF} it follows that
\begin{equation}\label{L3}
L_3\leq \frac16 \int_0^t \tau^{2-s}\|\sqrt{\rho}D_t^2 u\|_{L^2}^2d\tau+c\left(\|u_0\|_{H^{1+s}},T\right).
\end{equation}
Noticing that $\nabla\cdot D^2_t u=\nabla \cdot \left(u\cdot\nabla u_t +D_t u\cdot\nabla u+u\cdot D_t\nabla u\right)$, integration by parts twice in $L_2$ yields the following
\begin{equation}\label{L2aux}
\begin{aligned}
L_2&=\!\!\int_0^t\! \!\tau^{2-s}\!\!\!\int_{\R^2} \!\left(\nabla\!\cdot\! D_t^2 u\right)\! D_t p\hspace{0.05cm} dx d\tau=\!\!\!\int_0^t\!\! \tau^{2-s}\!\!\int_{\R^2} \!\nabla\!\cdot\! \left(u\cdot\nabla u_t\! +\!D_t u\cdot\nabla u\!+\!u\cdot D_t\nabla u\right)\! D_t p\hspace{0.05cm} dx d\tau\\
&=-\int_0^t \tau^{2-s}\int_{\R^2}\left(u\cdot\nabla u_t +D_t u\cdot\nabla u+u\cdot D_t\nabla u\right)\cdot \nabla D_t p  \hspace{0.05cm}dxd\tau\\
&\leq \int_0^t \tau^{2-s}\|u\cdot\nabla u_t +D_t u\cdot\nabla u+u\cdot D_t\nabla u\|_{L^2}\|\nabla D_t p\|_{L^2}d\tau.
\end{aligned}
\end{equation}

We now use the equation \eqref{Dt2} to estimate the high-order term $\nabla D_tp$. First, we notice that
\begin{equation}\label{Larayaux}
\mathbb{P}\Delta D_t u=\Delta D_t u-\nabla\nabla\cdot D_t u,
\end{equation}
so that the equation rewrites as
\begin{equation*}
-\mathbb{P}\Delta D_tu +\nabla D_tp=\nabla\nabla\cdot D_t u-\rho D_t^2 u+F.
\end{equation*}
Therefore we obtain 
\begin{equation*}
\|\mathbb{P}\Delta D_t u\|_{L^2}+\|\nabla D_t p\|_{L^2}\leq c\|\sqrt{\rho}D_t^2 u\|_{L^2}+\|F\|_{L^2}+\|\nabla \nabla \cdot D_t u\|_{L^2},
\end{equation*}
and using \eqref{Larayaux} we write
\begin{equation}\label{high}
\|\Delta D_t u\|_{L^2}+\|\nabla D_t p\|_{L^2}\leq c\|\sqrt{\rho}D_t^2 u\|_{L^2}+\|F\|_{L^2}+2\|\nabla \nabla \cdot D_t u\|_{L^2}.
\end{equation}
If we denote 
\begin{equation*}
G=u\cdot \nabla u_t+D_tu\cdot\nabla u+u\cdot D_t\nabla u,
\end{equation*}
going back to \eqref{L2aux}, estimate \eqref{high} provides the following
\begin{equation*}
L_2 \leq c\int_0^t \tau^{2-s}\|\sqrt{\rho}D_t^2 u\|_{L^2}\|G\|_{L^2}d\tau+\int_0^t \tau^{2-s}\left(\|F\|_{L^2}+2\|\nabla\nabla\cdot D_t u\|_{L^2}\right)\|G\|_{L^2}d\tau.
\end{equation*}
By Young's inequality it is possible to obtain
\begin{equation}\label{L2ini}
L_2\leq \frac16\int_0^t \tau^{2-s}\|\sqrt{\rho}D_t^2 u\|_{L^2}^2+c\int_0^t \tau^{2-s}\left(\|F\|_{L^2}^2+\|\nabla\nabla\cdot D_t u\|_{L^2}^2+\|G\|_{L^2}^2\right)d\tau.
\end{equation}
The incompressibility condition yields
\begin{equation*}
\|\nabla\nabla\cdot D_t u\|_{L^2}^2=\|\nabla (\nabla u\cdot\nabla u^*)\|_{L^2}^2\leq c \|\nabla u \cdot\nabla^2 u\|_{L^2}^2\leq c\|u\|_{H^{2+\eps}}^2\|\nabla^2 u\|_{L^2}^2.
\end{equation*}
Hence from \eqref{apriori1} and \eqref{uH2mu} we find that
\begin{equation}\label{tri1}
\int_0^t \tau^{1-s}\|\nabla\nabla\cdot D_t u\|_{L^2}^2\leq c\left(\|u_0\|_{H^{1+s}},T\right).
\end{equation}
On the other hand, the bound \eqref{uLinf} allows us to write
\begin{equation*}
t^{1-s}\|G\|_{L^2}^2\leq c \hspace{0.05cm}t^{1-s}\left(\|\nabla u_t\|_{L^2}^2+\|D_t u\|_{L^2}^2\|u\|_{H^{2+\eps}}^2+\|\nabla^2 u\|_{L^2}^2\right).
\end{equation*}
which joined to \eqref{Dtuestimates} and \eqref{apriori1} yields
\begin{equation*}
t^{1-s}\|G\|_{L^2}^2\leq c \left(\hspace{0.05cm}t^{1-s}\|\nabla u_t\|_{L^2}^2+\|u\|_{H^{2+\eps}}^2+1\right),
\end{equation*}
so  we conclude using again \eqref{apriori1} and \eqref{uH2mu} that
\begin{equation}\label{tri2}
\int_0^t \tau^{1-s}\|G\|_{L^2}^2d\tau\leq c \left(\|u_0\|_{H^{1+s}},T\right).
\end{equation}
If we introduce the bounds \eqref{boundF}, \eqref{tri1} and \eqref{tri2} in \eqref{L2ini} we get that
\begin{equation}\label{L24}
L_2\leq \frac16\int_0^t \tau^{2-s}\|\sqrt{\rho}D_t^2 u\|_{L^2}^2+c \left(\|u_0\|_{H^{1+s}},T\right).
\end{equation}
Finally, the term $L_4$ is bounded by
\begin{equation*}
L_4\leq \int_0^t \tau^{2-s}\|\nabla D_t u\|_{L^2}\|u\|_{H^{1+\eps}}\|\nabla^2 D_t u\|_{L^2}d\tau\leq c\int_0^t \tau^{2-s}\|\nabla D_t u\|_{L^2} \|\nabla^2 D_t u\|_{L^2}d\tau,
\end{equation*}
taking into account \eqref{uLinf}. Estimate \eqref{high}  gives
\begin{equation*}
L_4 \leq c\int_0^t \tau^{2-s}\|\nabla D_t u\|_{L^2}\|\sqrt{\rho}D_t^2 u\|_{L^2}d\tau\!+\!c\int_0^t\! \tau^{2-s}  \|\nabla D_t u\|_{L^2}\left(\|F\|_{L^2}\!+\!2\|\nabla\nabla\cdot D_t u\|_{L^2}\right)d\tau.
\end{equation*}
As in the bound of $L_2$, by Young's inequality we have that 
\begin{equation*}
L_4\leq \frac16\int_0^t \tau^{2-s}\|\sqrt{\rho}D_t^2 u\|_{L^2}^2+c\int_0^t \tau^{2-s}\|\nabla D_t u\|_{L^2}^2+c\left(\|u_0\|_{H^{1+s}},T\right),
\end{equation*}
and \eqref{Dtuestimates} implies
\begin{equation}\label{L42}
L_4\leq \frac16\int_0^t \tau^{2-s}\|\sqrt{\rho}D_t^2 u\|_{L^2}^2+c\left(\|u_0\|_{H^{1+s}},T\right).
\end{equation}

Introducing the bounds \eqref{L1}, \eqref{L3}, \eqref{L24} and \eqref{L42} in \eqref{Dt2balance}, we conclude that
\begin{equation*}
t^{2-s}\|\nabla D_t u\|_{L^2}^2+\int_0^t \tau^{2-s}\|\sqrt{\rho}D_t^2 u\|_{L^2}^2d\tau\leq c\left(\|u_0\|_{H^{1+s}},T\right).
\end{equation*}
Recalling \eqref{high}, \eqref{boundF} and \eqref{tri1} we find in addition that
\begin{equation*}
\int_0^t \tau^{2-s}\|D_t u\|_{H^2}^2d\tau\leq c\left(\|u_0\|_{H^{1+s}},T\right).
\end{equation*}
By Sobolev interpolation we note that in particular we have for $\tilde{s}\in(0,s)$,
\begin{equation}\label{utLpH1s}
D_t u \in L^p(0,T; H^{1+\tilde{s}}), \hspace{0.5cm}1\leq p<\frac{2}{2-(s-\tilde{s})}.
\end{equation}

\vspace{0.1cm}
\subsection{Critical regularity for $u$}
\vspace{0.3cm}

In this section we will conclude that $\displaystyle u\in L^p(0,T; W^{2,\infty})$.
 From \eqref{laplace} we have that
\begin{equation*}
\nabla^2 u=\nabla^2 \Delta^{-1}\mathbb{P}\left(\rho D_t u\right),
\end{equation*}
where $\mathbb{P}f_i= f_i-R_i R_j f_j$. The operators $\nabla^2\Delta^{-1}\mathbb{P}$ are Fourier multipliers and therefore can be written as convolutions
\begin{equation}\label{conv}
\partial_k\partial_l\Delta^{-1}\mathbb{P}f_i(x)=\left(K_{klij}\star f_j\right)(x),
\end{equation}
with kernels given by
\begin{equation}\label{kernelfou}
K_{klij}(x)=\mathcal{F}^{-1}\left(\frac{\xi_k\xi_l}{|\xi|^2}\left(\delta_{ij}-\frac{\xi_i\xi_j}{|\xi|^2}\right)\right)(x).
\end{equation}
By symmetries it suffices to consider the following three cases:
\begin{equation*}
\begin{aligned}
\partial_1^2\Delta^{-1}\mathbb{P}f_1(x)&=\left(\tilde{K}_{111j}\star f_j\right)(x)+\frac18f_1(x),\\
\partial_1^2\Delta^{-1}\mathbb{P}f_2(x)&=\left(\tilde{K}_{112j}\star f_j\right)(x)+\frac38 f_2(x),\\
\partial_1\partial_2 \Delta^{-1}\mathbb{P}f_1 (x)&=\left(\tilde{K}_{121j}\star f_j\right)(x)-\frac18f_2(x).
\end{aligned}
\end{equation*}
where the kernels $\tilde{K}_{klij}$ are even and have zero mean on circles. They can be computed explicitly as they correspond to sums of second and fourth order Riesz transforms (see Chapter 3.3 in \cite{STEIN}).
For simplicity we will denote by $K$ the kernels $\tilde{K}_{klij}$ and we rewrite the above equations as singular integral operators plus identities as follows
\begin{equation}\label{hol}
\nabla^2 u= \textup{SIO}\left( \rho D_t u \right)+c\hspace{0.02cm}\rho D_t u.
\end{equation}
Then we decompose as follows
\begin{equation*}
\begin{aligned}
\textup{SIO}\left( \rho D_t u \right)&=\rho_2 \hspace{0.05cm} \textup{SIO}\left(D_t u \right)\\
&\quad+ \left(\rho_1-\rho_2\right)\int_{D(t)}K(x-y)\cdot\left(D_tu (y,t)-D_tu(x,t)\right)dy\\
&\quad+\left(\rho_1-\rho_2\right)\textup{SIO}\left(1_{D(t)}\right) D_tu(x,t)=M_1+M_2+M_3.
\end{aligned}
\end{equation*}
By Sobolev embedding and \eqref{utLpH1s} we get that
\begin{equation}\label{DtuHolder}
D_t u \in L^p(0,T;C^{\tilde{s}}),
\end{equation} 
hence we deduce that
\begin{equation*}
M_1\leq c\left( \|D_tu\|_{C^{\tilde{s}}}(t)+\|D_tu\|_{L^2}(t) \right),
\end{equation*}
and analogously
\begin{equation*}
M_2\leq c\|D_tu\|_{C^{\tilde{s}}}(t).
\end{equation*}
As by Theorem \ref{Case1} $\rho(t)$ is a $C^{1+\gamma}$ patch for all $t\geq0$, $\gamma\in(0,1)$, and the fact that the kernels in the singular integral operators are even, it is possible to obtain (see \cite{CONST} for more details)
\begin{equation*}
M_3\leq c\|D_tu\|_{L^\infty}(t).
\end{equation*}
We therefore conclude that
\begin{equation}\label{siobound}
\|\textup{SIO}\left( \rho D_t u \right)\|_{L^p_T(L^\infty)}\leq  c\left(\|u_0\|_{H^{1+s}},T\right),
\end{equation}
and hence \eqref{hol} gives
\begin{equation*}
u\in L^p(0,T; W^{2,\infty}),\hspace{0.5cm}1\leq p<\frac{2}{2-(s-\tilde{s})}.
\end{equation*}
\vspace{0.3cm}

\begin{rem}
We have obtained $u\in L^2(0,T; H^{2+\mu})$ with $\mu<\min\{s,1/2\}$. In the case $s<1/2$ we can also get $u\in L^q(0,T; H^{2+\delta})$ with $\delta\in(s,1/2)$ and $q\in[1,2/(1+\delta-s))\subset[1,2)$. This is achieved by rewriting  the equation as in \eqref{laplace} to take advantage of the smoothing properties of the Laplace equation. If $s<1/2$, from  $D_t u\in L^1(0,T;H^{1+\tilde{s}})\cap L^2(0,T; H^{\tilde{s}})$ one finds by  interpolation that
$D_t u \in L^q(0,T;H^{\delta})$. Since $\rho$ is a multiplier in $H^\delta$ for any $\delta\in(s,1/2)$, from standard properties of the Laplace equation it follows that
\begin{equation*}
u \in L^q(0,T;H^{2+\delta}).
\end{equation*}	
\end{rem}

\qed

\section{Persistence of $C^{2+\gamma}$ regularity}\label{sec:4}

This section is devoted to show that $C^{2+\gamma}$ regularity is preserved globally in time.

\begin{thm}\label{Case3}
	Assume $\gamma\in(0,1)$, $s\in(0,1-\gamma)$, $\rho_1, \rho_2>0$. Let $D_0\subset \R^2$ be a bounded simply connected domain with boundary $\partial D_0\in C^{2+\gamma}$, $u_0\in H^{1+\gamma+s}$ a divergence-free vector field and  
	\begin{equation*}
	\rho_0(x)=\rho_1 1_{D_0}(x)+\rho_2 1_{D_0^c}(x).
	\end{equation*}
	Then, there exists a unique global solution $(u,\rho)$ of \eqref{INH}-\eqref{NavierStokes} such that
	$$\rho(x,t)=\rho_1 1_{D(t)}(x)+\rho_2 1_{D(t)^c}(x) \hspace{0.2cm}{\rm{and}} \hspace{0.2cm} \partial D\in C([0,T];C^{2+\gamma}).$$
	Moreover, 
	\begin{equation*}
	\begin{aligned}
	u&\in C([0,T];H^{1+\gamma+s})\cap L^2(0,T;H^{2+\mu})\cap L^p(0,T;W^{2,\infty}),\\
	t^{\frac{1-(\gamma+s)}{2}} u_t&\in L^\infty(\left[0,T\right];L^2)\cap L^2(\left[0,T\right];H^1),\\
	t^{\frac{2-(\gamma+s)}{2}} D_tu&\in L^\infty(\left[0,T\right];H^1)\cap L^2(\left[0,T\right];H^2),
	\end{aligned}
	\end{equation*}
	for any $T>0$, $\mu<\min\{1/2,\gamma+s
	\}$ and $p\in[1,2/(2-(\gamma+s)))$. If $\gamma+s<1/2$ it also holds that $u\in L^q(0,T;H^{2+\delta})$ for any $\delta \in(\gamma+s,1/2)$ with $q\in[1,2/(1+\delta-\gamma-s))$.
\end{thm}

Proof: Since $\gamma+s\in(0,1)$, the estimates on $u$, $u_t$ and $D_tu$ follow as in Theorem \ref{Case2}. 
We now describe the patch using a level-set function $\varphi$: 
\begin{equation*}
\partial_t\varphi+u\cdot\grad\varphi=0,\hspace{0.7cm}
\varphi(x,0)=\varphi_0(x),
\end{equation*}
\begin{equation*}
D_0=\{x\in\R^2: \varphi_0(x)>0\},
\end{equation*}
so that at time $t$, $D(t)=X(t, D_0)=\{x\in\R^2: \varphi(x,t)>0\}$.
Then, the vector field given by $W(t)=\grad^\perp \varphi(t)$ is tangent to the patch and evolves as follows
\begin{equation}\label{tangent}
\begin{aligned}
\partial_t W +u\cdot \grad W=W\cdot\grad u,\hspace{0.5cm}
W(0)=\grad^\perp \varphi_0.
\end{aligned}
\end{equation}
To control $C^{2+\gamma}$ regularity of $\partial D(t)$ we shall ensure that $\grad W$ remains in $C^{\gamma}$.  By differentiating \eqref{tangent} one obtains

\begin{equation*}
\partial_t \grad W + u\cdot\grad(\grad W) =W\cdot\nabla^2 u + \grad W\cdot \grad u +\grad u \cdot \grad W.
\end{equation*}
Since $u$ is Lipschitz, the following estimate holds for all $t\in[0,T]$:
\begin{equation*}
\begin{aligned}
\|\grad W\|_{C^\gamma}(t)&\!\leq\!\! \|\grad W_0\|_{C^\gamma}e^{c\int_0^t\|\grad u\|_{L^\infty}d\tau}\!+\!e^{c\int_0^t\|\grad u\|_{L^\infty}d\tau}\!\!\int_0^t\! \left( \|W\cdot\grad^2u\|_{C^\gamma}\!+\!2\|\grad W\|_{L^\infty} \|\grad u\|_{C^\gamma}\!\right)d\tau.
\end{aligned}
\end{equation*}
From this and previous estimates we get that
\begin{equation*}
\|\grad W\|_{C^\gamma}(t)\leq c_1(T)+c_2(T)\int_0^t \|W\cdot\grad^2u\|_{C^\gamma}(\tau)d\tau.
\end{equation*}
Therefore the result is obtained once we prove that $W\cdot\nabla^2 u\in L^1(0,T; C^\gamma)$.
\\

Applying Fourier transform in \eqref{conv} gives that 
\begin{equation*}
\mathcal{F}\left(W_k\partial_k\partial_l u_i\right)(\xi)=\hat{W}_k(\xi)\star \mathcal{F}\left(\partial_k\partial_l u_i\right)(\xi)=\hat{W}_k(\xi)\star\left[ \frac{\xi_k\xi_l}{|\xi|^2}\left( \delta_{ij}-\frac{\xi_i\xi_j}{|\xi|^2} \right)\mathcal{F}\left(\rho D_tu_j\right)(\xi) \right].
\end{equation*}
Using the notation \eqref{kernelfou}, we introduce the following splitting
\begin{equation*}
\begin{aligned}
\mathcal{F}\left(W_k\partial_k\partial_l u_i\right)(\xi)&=\left( \hat{K}_{ijkl}\mathcal{F}\left(\rho D_t u_j\right) \right)\star \hat{W}_k\hspace{0.02cm}(\xi)-\hat{K}_{ijkl}(\xi)\mathcal{F}\left(W_k\rho D_t u_j\right)(\xi)\\
&\quad+\hat{K}_{ijkl}(\xi)\mathcal{F}\left(W_k \rho D_t u_j\right)(\xi)-\left(\hat{K}_{ijkl}\mathcal{F}\left(\rho W_k\right)\right)\star \mathcal{F}\left(D_tu_j\right)(\xi)\\
&\quad+\left(\hat{K}_{ijkl}\mathcal{F}\left(\rho W_k\right)\right) \star\mathcal{F}\left(D_tu_j\right)(\xi).
\end{aligned}
\end{equation*}
We note that since $W$ is tangent to the density patch, the last term vanish 
\begin{equation*}
\hat{K}_{ijkl}\mathcal{F}\left(\rho W_k\right)(\xi)=-i\left(\frac{\xi_l}{|\xi|^2}\left(\delta_{ij}-\frac{\xi_i\xi_j}{|\xi|^2}\right)\right) \mathcal{F}\left(\partial_k \left(W_k\rho\right)\right)(\xi)=0.
\end{equation*}
Hence the previous splitting writes as
\begin{equation*}
W(x,t)\cdot \nabla^2 u(x,t)=I_1+I_2,
\end{equation*}
where
\begin{equation*}
\begin{aligned}
I_1&=\int_{\R^2}K(x-y)\cdot\left(W(x,t)-W(y,t)\right)\rho(y,t)D_t u(y,t)dy,\\
I_2&=\int_{\R^2} K(x-y)\cdot W(y,t)\rho(y,t)\left(D_t u(y,t)-D_t u(x,t)\right)dy.
\end{aligned}
\end{equation*}
The Lemma in Appendix of \cite{CONST} yields the following
\begin{equation*}
\begin{aligned}
\|I_1\|_{C^{\gamma}}&\leq c \|W\|_{C^{\gamma}}(t) \left(\|\rho D_t u\|_{L^\infty}(t)+\|\textup{SIO}\left( \rho D_t u\right)\|_{L^\infty}(t)\right),\\
\|I_2\|_{C^{\gamma}}&\leq c \|D_t u\|_{C^{\gamma}}(t) \left(\|\rho W\|_{L^\infty}(t)+\|\textup{SIO}\left( \rho W \right)\|_{L^\infty}(t)\right).
\end{aligned}
\end{equation*}
Proceeding as in  \eqref{siobound} it is possible to find that 
\begin{equation*}
\begin{aligned}
\|I_1\|_{L^1_T(C^{\gamma})}&\leq c \left( \|u_0\|_{H^{1+\gamma+s}},T\right),\\
\|I_2\|_{L^1_T(C^{\gamma})}&\leq c\left( \|u_0\|_{H^{1+\gamma+s}},T\right).
\end{aligned}
\end{equation*}
From the above estimates we finally conclude that
\begin{equation*}
\|W\cdot \nabla^2 u\|_{L^1_T(C^{\gamma})}\leq c\left( \|u_0\|_{H^{1+\gamma+s}},T\right).
\end{equation*} 

\qed

\subsection*{{\bf Acknowledgments}}
This research was partially supported by the project P12-FQM-2466 of Junta de Andaluc\'ia, Spain, 
 the grant MTM2014-59488-P (Spain) and by the ERC through the Starting Grant project H2020-EU.1.1.-639227. 
EGJ was supported by MECD FPU grant from the Spanish Government.


\vspace{1cm}
\newpage

\begin{quote}
	\begin{tabular}{ll}
		\textbf{Francisco Gancedo}\\
		{\small Departamento de An\'{a}lisis Matem\'{a}tico $\&$ IMUS}\\
		{\small Universidad de Sevilla}\\
		{\small C/ Tarfia s/n, Campus Reina Mercedes,}\\
		{\small	41012 Sevilla, Spain}\\
		{\small Email: fgancedo@us.es}
	\end{tabular}
\end{quote}

\begin{quote}
	\begin{tabular}{ll}
		\textbf{Eduardo Garc\'ia-Ju\'arez}\\
		{\small Departamento de An\'{a}lisis Matem\'{a}tico $\&$ IMUS}\\
		{\small Universidad de Sevilla}\\
		{\small C/ Tarfia s/n, Campus Reina Mercedes,}\\
		{\small	41012 Sevilla, Spain}\\
		{\small Email: eduardogarcia@us.es}
	\end{tabular}
\end{quote}

\end{document}